\documentclass{article}
\usepackage{graphicx} % Required for inserting images

\usepackage[a4paper, total={5.8in, 8.5in}]{geometry}
\usepackage{amsmath,amsthm,bm}
\usepackage{amsfonts}
\usepackage{mathtools}
\usepackage{ stmaryrd }
\usepackage{tikz-cd}
\usepackage{ upgreek }
\usepackage{ amssymb }
\usepackage{dsfont}
\usepackage{hyperref}
\usepackage{yfonts}
\usepackage{hyperref}
\usepackage{comment}
\usepackage[colorinlistoftodos]{todonotes}
\usepackage{bbm}

\usepackage{thmtools}
\usepackage{cleveref}
\usepackage{nameref}
\Crefname{rem}{Remark}{Remarks}
\Crefname{alphthm}{Theorem}{Theorems}
\Crefname{thm}{Theorem}{Theorems}
\Crefname{equation}{Equation}{Equations}

\usepackage{titling}
\usepackage[affil-it]{authblk}
\usepackage{orcidlink}
\newcommand{\email}[1]{\href{mailto:#1}{\texttt{#1}}}

\usepackage[
backend=biber,
maxnames=10,
style=numeric,
sorting=anyt
]{biblatex}
\makeatletter
\newcommand{\ostar}{\mathbin{\mathpalette\make@circled\star}}
\newcommand{\make@circled}[2]{%
  \ooalign{$\m@th#1\smallbigcirc{#1}$\cr\hidewidth$\m@th#1#2$\hidewidth\cr}%
}
\newcommand{\smallbigcirc}[1]{%
  \vcenter{\hbox{\scalebox{0.77778}{$\m@th#1\bigcirc$}}}%
}
\makeatother

\makeatletter
\DeclareRobustCommand{\nsext}[1]{%
  \binrel@{#1}% compute the type
  \binrel@@{%
    {\vphantom{#1}}^*% the asterisk at the proper height
    \kern-\scriptspace % remove the script space
    \csname mkern@\detokenize{#1}\endcsname % additional kerning
    {#1}% the symbol
  }%
}
\newcommand{\defineextkern}[2]{%
  \@namedef{mkern@\detokenize{#1}}{\mkern#2}%
}
\makeatother

\newcommand\restr[2]{{% we make the whole thing an ordinary symbol
  \left.\kern-\nulldelimiterspace % automatically resize the bar with \right
  #1 % the function
  \littletaller % pretend it's a little taller at normal size
  \right|_{#2} % this is the delimiter
  }}

\newcommand{\littletaller}{\mathchoice{\vphantom{\big|}}{}{}{}}

\newcommand{\nf}[1][]{%
  \mathrel{
    \mathop{
      \vcenter{
        \hbox{\oalign{\noalign{\kern-.3ex}\hfil$\vert$\hfil\cr
              \noalign{\kern-.7ex}
              $\smile$\cr\noalign{\kern-.3ex}}}
      }
    }\displaylimits_{#1}
  }
}

\usepackage{afterpage}

\theoremstyle{definition}
\newtheorem{thm}{Theorem}[section]
\newtheorem{lemma}[thm]{Lemma}

\newtheorem{rmrk}[thm]{Remark}

\newtheorem{coroll}[thm]{Corollary}

\newtheorem{fact}[thm]{Fact}
\newtheorem{defn}[thm]{Definition}

\newtheorem*{notation}{Notation}

\newtheorem{question}[thm]{Question}
\newtheorem*{claim}{Claim}
\newtheorem*{ass}{Assumptions}

\newtheorem{alphthm}{Theorem}

\newtheoremstyle{named}{}{}{}{}{\bfseries}{.}{.5em}{\thmnote{#3}}
\theoremstyle{named}

\let\oldqed\qedsymbol

\newcommand{\qedclaim}{\mbox{$\underset{\textsc{claim}}{\oldqed}$}}
\makeatletter
\newenvironment{claimproof}[1][\it Proof of Claim]{
\let\qedsymbol\qedclaim
  \par
  \pushQED{\qed}%
  \normalfont \topsep6\p@\@plus6\p@\relax
  \trivlist
\item[\hskip\labelsep
  \upshape
  #1\@addpunct{.}]\ignorespaces
}{%
  \popQED\endtrivlist\@endpefalse
}
\makeatother
\let\qedsymbol\oldqed

\newcommand{\tp}{\text{tp}}

\makeatletter
\newcommand{\subjclass}[2][2020]{%
  \let\@oldtitle\@title%
  \gdef\@title{\@oldtitle\footnotetext{\hspace*{-2em}#1 \emph{Mathematics subject classification.} #2}}%
}
\newcommand{\keywords}[1]{%
  \let\@@oldtitle\@title%
  \gdef\@title{\@@oldtitle\footnotetext{\hspace*{-2em}\emph{Keywords:} #1.}}%
}

\makeatother

\def\R{\mathbb{R}}
\def\Q{\mathbb{Q}}
\def\N{\mathbb{N}}
\def\Z{\mathbb{Z}}

\def\U{\mathbb{U}}

\author[1]{Alessandro Vegnuti\orcidlink {0009-0005-2490-720X}%
  \thanks{email: \email{alessandro.vegnuti@unimi.it}}}

\affil[1]{\small Dipartimento di Matematica, Università  di Milano, Via Saldini 50, 20133 Milano, Italy and Mathematisches Institut, Albert-Ludwigs-Universität Freiburg, D-79104 Freiburg,
Germany}

\title{Solving polynomials in large pseudo-finite sets}

\date{\today}



\begin{document}

\maketitle

\begin{abstract}
    Using techniques from model theory, we prove analogs of Rado's and S{\'a}rk{\"o}zy's theorems on solutions to polynomial equations in dense sets of definable groups and fields equipped with well-behaved measures. As a special case, we adapt our results to pseudo-finite structures, to deduce non-qualitative asymptotic results. 
\end{abstract}

\section{Introduction}

One of the main problems in additive combinatorics is to determine which kind of structures always appear in large enough sets. For example, in 1975 Szemerédi \cite{szemeredi1975sets} answered positively a conjecture of Erdős and Turán, proving that any large subset of the integers contains arbitrarily long arithmetic progressions, extending Roth's theorem \cite{roth1953certain} on 3-term arithmetic progressions.

Observe that if $(a_1, a_2, a_3)$ is a tuple in $X^3$ that solves $x + y = 2z$ and we pose $g=a_3 - a_1$, then $\{a_1, a_1 + g, a_1 + 2g\}$ is a 3-term arithmetic progression contained in $X$. Hence, to find such a progression, it suffices to solve a linear equation (or a system, in the case of more terms) within a fixed set. This requires a deep study of sets of solutions to polynomials in several variables.

For the linear case, Rado proved in \cite{Rado} that a linear polynomial is \textit{partition regular}, i.e. it admits a monochromatic solution given any finite coloring of $\N$, if and only if a sum of its coefficients is zero (we call this \textit{Rado's condition}). However, a density version of this fact has not yet been written down (though the recent paper of Hernandez \cite{hernandez2025infinite} considers linear patterns in translates of dense sets). For the nonlinear case, the most notable result is due to S{\'a}rk{\"o}zy \cite{sarkozy1978difference}, who proved that in any dense subset of the integers we can find two numbers whose difference is a perfect square.  

Notice that all these problems can be re-formulated in a finite setting: if $N$ is a positive integer and we take $X$ a subset of $\{1,\dots , N\}$ such that its cardinality is large enough, then we expect to find some structure in $X$. Many results have been obtained in this direction, focused both on the qualitative (which patterns have such a property) and quantitative (how small can a subset be to avoid patterns) sides of the questions: among others, we mention here \cite{RadoFinite} for Rado's theorem over finite fields, and \cite{mccutcheon2003sarkozy} for S{\'a}rk{\"o}zy's theorem, again over finite fields.

A new direction was developed in recent years: starting from the striking result of Hrushovski \cite{hrushovski2012stable}, which led to the characterization of approximate subgroups of local groups \cite{breuillard2012structure}, several applications of model theory to combinatorics have been found (see e.g. \cite{martin2026supersimplicity,massicot2015approximate,palacin2020compactifications} and \cite{terry2019stable}). 

In this paper, we further explore these connections, focusing on solutions of polynomials in large sets. Our first result involves solving linear equations in translates of sets which are dense with respect to a well-behaved measure, giving a model-theoretical analogue of Rado's theorem, and generalizes a result of Martin-Pizarro and Palacín \cite[Remark 3.15]{martin2024complete}. Here $G$ is a definable group in some saturated structure $\mathbb U$, using parameters from an elementary substructure  $ M \prec \mathbb U$; furthermore, we assume that in every Cartesian power $G^n$ we have an additive, symmetric measure $\mu_n$ which is invariant under translation and such that $\mu_n(B)>0$ if and only if $\mu_n(m\cdot B)>0$ for every non-zero integer $m$. We also assume that the family $\{\mu_n\}_{n\in \N}$ satisfies a form of Fubini (cf. Section 2.2).

\begin{alphthm}[\Cref{teo rado su traslati}]
    Let $L=n_1x_1 + \dots + n_s x_s$ be a linear polynomial without constant term with a sum of its coefficients equals to 0. 
    
    If $X$ is an $M$-definable subset of $G$ of positive measure, then for some $m$ in $G(M)$ the set \begin{equation*}
        \{(x_1,\dots, x_{s-1})\in (X+m)^{s-1} \mid \exists x_s\in (X+m): L(x_1,\dots, x_s)=0\}
    \end{equation*}
    has positive $\mu_{s-1}$-measure.
\end{alphthm}

The main model-theoretical ingredients will be a model-theoretical version of Croot-Sisask’s almost-periodicity first proved in \cite{martin2024complete}, and Hrushovski's stabilizer theorem, linking the measure-theoretical stabilizer of types to the connected component of the group.

We then move on to the study of fields. Here we assume the measures to be not only invariant under additive, but also multiplicative shifts. We generalize a result of Pillay, Scanlon and Wagner \cite[Lemma 3.2]{pillay1998supersimple} proving that the additive and multiplicative components intersect largely, allowing us to use the stabilizer theorem for solving nonlinear equations. For example, we deduce the following. 

\begin{alphthm}[\Cref{4:co:EquazioniInGenerale}]
    Let $X$ be a definable subset of an $M$-definable field $F$ such that $\mu(X)>0$. For all elements $c_1, \dots, c_l, d_1,\dots$ and  $d_{s} $ of $F$ such that $\sum_{i=1}^l c_i=0$, there exists some $g$ in $F(M)$ such that the set \begin{equation*}
        \{(x_1, \dots , x_{l-1}, y_1,\dots, y_s )\in (X + g)^{ l+ s- 1} \mid \exists x_{l}\in (X + g) :  \sum_{i=1}^l c_i\cdot x_i^n + \sum_{j=1}^s d_j\cdot y_j^m =0\}\end{equation*} 
    has positive $\mu_{l+ s- 1}$ measure.
\end{alphthm}

As for a natural application, we consider pseudo-finite groups and fields: these are ultraproducts of finite groups and fields, and their study allows for an understanding of the asymptotic behavior of finite structures. In particular, a standard argument produces the following asymptotic, qualitative analogs of the the theorems above.

\begin{alphthm}[\Cref{3:co:RadoPerFiniti} and  \Cref{4:co:NonLinearSuiFiniti}]
    Fix natural numbers $n $ and $m $ as well as a real positive number $\epsilon$. For every linear polynomial $L(x_1,\dots, x_s)=n_1x_1 + \dots + n_s x_s$ over $\Z$ with $\sum_{i=1}^{l-1}n_i=-n_l$, there exist values $c$ and $\delta$ such that the following holds.

    \begin{itemize}
        \item If $G$ is a finite group with $|G|>c$ coprime with $n_l\cdot n_s\cdot (\sum_{j=l+1}^s n_j)$, for every subset $X$ of $G$ of cardinality at least $\epsilon c$ there exists $g$ such that the set \begin{equation*}
        \{(x_1,\dots, x_{s-1})\in (X+g)^{s-1} \mid \exists x_s \in (X + g) : L(x_1,\dots, x_s)=0\}
    \end{equation*}
    has cardinality at least $\delta c^{s-1}$.

    \item For every finite field $K$ with $|G|>c$ and for every subset $X$ of $K$ with cardinality at least $\epsilon c$, as well as elements $c_1,\dots, c_l, d_{1},\dots, d_{s-1}$ and $d_s$ with $\sum_{i=1}^l c_i =0$, the set \begin{equation*}
        \{(x_1, \dots , x_{l-1}, y_1,\dots, y_{s})\in (X+ g)^{l +s -1} \mid \exists x_{l} \in (X + g):  \sum_{i=1}^{l} c_i\cdot x_i^n + \sum_{j=1}^s d_j\cdot y_j^m =0\}\end{equation*} has cardinality at least $\delta c^{l +s -1}$.
    \end{itemize}

\end{alphthm}

\paragraph{Acknowledgments}
The author would like to thanks his supervisor Amador Martin-Pizarro for his invaluable support and guidance during the preparation of these notes.
\paragraph{AI disclousure}
AI was only used for some formatting of LaTeX. In particular, all the mathematics herein is human-generated.

\section{Model-theoretical preliminaries: wide types, connected components and stabilizers}

In this section, we recall the main notions we will need later in the paper. We will assume some familiarity with basic concepts of model theory such as types and saturation; for a more detailed exposition, we suggest \cite{tent2012course}, among others. To fix the terminology, we work in an infinite, sufficiently saturated model $\mathbb U$ of some complete first-order theory in a countable language, and $M$ is a countable elementary submodel. We use $A, B$ or $C$ for small sets of parameters, and $X, Y$ or $Z$ for definable sets. 

\begin{comment}
We work in a complete theory in a countable language $\mathcal L$, with $\mathbb U$ as an infinite, sufficiently saturated model and $M$ an elementary submodel. If not specified, a formula will have parameters coming from $\mathbb U$. Notice that we will make no distinction between a formula and the set it defines. Moreover, all parameters and variables are intended to be tuples, unless explicitly stated.

A partial type over $B\subset \mathbb U$ is a family of formulas with parameters coming from $B$ which is finitely consistent. A (complete) type is a partial type which is maximal by inclusion, and any partial type can be completed. We denote with $S_n(B)$ the set of all complete types over $B$ in $n$ variables; if the number of variables is not essential, we omit it. We say that $c$ realizes the type $p$ if it satisfies all formulas in $p$; by saturation, every (partial) type over a small set of parameters $B$ is realized in $\mathbb U$. Viceversa, given a tuple $c$ we define its type $\tp(c/B)$ over $B$ to be the set of formulas $\varphi(x, b)$ with parameters in $B$ such that $\varphi(c, b)$ holds in $\mathbb U$.  

We endow $S_n(B)$ with a topology whose clopen sets are $$\{[\varphi(x, b)]\mid \varphi \text{ is an $\mathcal{L}$-formula, }b\in B\}$$ where $[\varphi(x, b)]=\{p\in S_n(B)\mid \varphi(x, b)\in p\}$. With this topology $S_n(B)$ is a Stone space, that is: a compact, 0-dimensional Hausdorff space.
\end{comment}

\subsection{Keisler measures}\label{subsec definability Keisler measures}

We recall that a Keisler measure in dimension $n$, as introduced in \cite{keisler1987measures}, is a finitely additive probability measure on the boolean algebra of definable subsets of $\mathbb U^n$. %To be more precise, we consider a family of measures $\{\mu_n\}_{n>0}$ for each dimension $n$, so the set defined by $\varphi(x, a)$ will be measured via $\mu_{|x|}$: however, we will usually refer to each one of them using $\mu$.
As a general reference, we invite the reader to consult Simon's book \cite{simon}.

\begin{defn}\label{def misura definibile}
    Let $\mu$ be a Keisler measure. We say that $\mu$ is $A$-definable if given any $\mathcal{L}$-formula $\varphi(x,y)$ with $|x|=n$ and any $\varepsilon >0$ there is a partition of the universe into formulas $\rho_1, \dots, \rho_m$ with parameters from $A$ such that if $b$ and $b'$ realize the same $\rho_i$, then $$\vert \mu(\varphi(x, b))- \mu(\varphi(x, b')) \vert \leq \varepsilon.$$
\end{defn}

It can be shown that $A$-definability of $\mu $ is equivalent to saying that $\mu $ is $A$-invariant and for every $\varphi(x,y)$ the map 
    \[
        \begin{array}{r@{\quad}c@{\quad}l}
f: & S_m(A) & \to\ \  \mathbb{R}\\
              & tp(b/A) & \mapsto \ \ \mu(\varphi(x,b))
\end{array}
\]
is continuous. 

\begin{rmrk}\label{rem b in Y}
    Suppose that the measure $\mu$ is definable over $A$ and let $\varphi(x,y)$ be an $\mathcal{L}$-formula such that $\mu(\varphi(x,b))=r$. For any real number $\varepsilon>0$ there exists a $A$-definable set $Y$ such that $$c\in Y \subset \{y \in \mathbb U \mid r - \varepsilon \leq  \mu(\varphi(x, y))\leq r +\varepsilon\}.$$
    Indeed, it suffices to consider the formula $\rho_i$ that contains $b$ (\Cref{def misura definibile}).
\end{rmrk}

Given a Keisler measure $\mu$ on $\mathbb U$, we can define a $\sigma$-additive Borel measure on $S(\mathbb U)$ such that, if $C=[\varphi(x)]$, then the measure of $C$ coincides with $\mu(\varphi(x))$. We denote this measure $\mu$ again. For a more detailed construction, we suggest \cite[Section 7]{simon}.

\subsection{Wide types}

Let $G$ be a definable group, that is, a definable subset of $\mathbb U^n$ which is a group and such that the group operation and inverse are both definable functions. We can assume that $G$ is definable over $M$ for $M$ a countable elementary submodel of $\mathbb U$. The interpretation of $G$ in $\mathbb U$ will be called just $G$, while in $M$ it will be specified $G(M)$. Consider a family $\{\mu_n\}_{n\in \N} $ of definable Keisler measures on $G^n$, invariant by translations and by taking the inverse. We also assume that they satisfy Fubini's condition: 
$$\mu_{|x| + |y|}(\varphi(x, y))=\int_{tp(b/M)\in S_{|b|}(M)}\mu_{|x|}(\varphi(x, b))d\mu(c)=\int_{tp(a/M)\in S_{|a|}(M)} \mu_{|y|}(\varphi(a, y))d\mu(a)$$
where the integral is considered with respect to the $\sigma$-additive measure on the space of types. If clear from the context, we will drop the subscript from $\mu_n$. As already mentioned, the map that sends $tp(b/M)$ to $\mu(\phi(x, b))$ is continuous, hence measurable: so, the above integral is well-defined. 

\begin{defn}
    A (partial) type over $A$ is wide if it does not contain any definable set of measure zero. A tuple $b$ is wide over $A$ if $\tp(b/A)$ is.
\end{defn}

\begin{fact}\label{F wr tuples}
The following holds:
\begin{itemize}
    \item Every wide partial type $\Sigma$ over $A$, closed under deductions, can be extended to a complete wide type over $B\supset A$ (by compactness applied to $\Sigma\cup \{\lnot \varphi(x, b)\mid \mu(\varphi(x, b))=0\ ,\ b\in B\}$);
    \item If $a$ is wide over $C$ and $b$ is wide over $C\cup \{b\}$, then $(a,b)$ is wide over $C$ (\cite[Lemma 1.10]{martin2024complete});
    \item The pair $(a, b)$ is wide over $C$ if and only if $(b, a)$ is (due to Fubini). 
\end{itemize}
    
\end{fact}

In general, it is not true that if $(a,b)$ wide over $C$ then $a$ is wide over $C\cup \{b\}$. %For example ....\todo{Candidate: to work like in DLO + operation; but about definability of the measure?}
Mainly for this reason, we need to introduce the stronger notion of random types. In order to do so, we need some preliminary definitions.

\begin{defn}
    Fix a set of parameters $C$. For every Borel set $Z\subset S_{n+m}(C)$ define $$Z(x, b )=\{p\in S_n(\mathbb U) \mid  \text{ for some } a\models p \upharpoonright_{C,b}, \ \tp(a,b/C)\in Z\}.$$
    Given $Z$ as above, we can define the function $$
    \begin{array}{r@{\quad}c@{\quad}l}
F^Z_{\mu, C}: & S_m(C) & \to\ \  \mathbb{R}\\
              & \tp(b/C) & \mapsto\ \  \mu(Z(x,b))
\end{array}
    $$ which is well-defined and continuous. 
\end{defn}

\begin{defn}
    Let $\mathcal B_M$ be the smallest Boolean algebra over $\bigcup_n S_n(M)$ that contains all the clopen sets and such that, for every $Z\in \mathcal B_M$:\begin{itemize}
        \item $\pi^{-1}(Z)\in \mathcal B_M$, where $\pi: S_m(M)\to S_n(M)$ ($n\leq m$) is the projection on some fixed variables;
        \item $(F_{\mu, M}^Z)^{-1}(\{0\})\in \mathcal B_M$.
    \end{itemize}
\end{defn}

Notice that all elements in $\mathcal B_M$ are Borel. 

\begin{defn}
    Let $B$ be a countable set of parameters. A tuple $a$ is random over $M\cup B$ if for every $Z(x,y)\in \mathcal B_M$ and every $b\in B^{\vert y \vert}$, if $\tp(a,b/M)\in Z$ then $\mu(Z(x,b))>0$. 
\end{defn}

An easy application of $\sigma$-additivity ensures that a Borel set of positive measure in the space of types must contain a random element. In particular, if $A$ is an $M\cup B$-definable set of positive measure, then it contains a random element over $M\cup  B$. Moreover, since all clopen sets are contained in $\mathcal B_M$, random types are in particular wide. As already mentioned, unlike wide types, random ones are symmetric, in the following sense.

\begin{fact}[{\cite[Lemma 1.15]{martin2024complete}}]
    The tuple $(a,b)$ is random over $M$ if and only if $a$ is random over $M$ and $b$ is random over $M\cup \{a\}$. 
\end{fact}

\begin{rmrk}
    Invariance by translations and by inverse of the measure ensures that if $a$ is random over $M\cup B$, then so are $a^{-1}$ and $a\cdot b$ (for $b$ in $ B$). 
    Moreover, an easy application of Fubini yields that a pair $(a,b)$ is random if and only if $(b,a)$ is.
\end{rmrk}
 Also, we can obtain symmetry for wide types, assuming that one of them is random.

\begin{fact}[{\cite[Lemma 1.16]{martin2024complete}}]\label{referee's lemma}
    Let $M$ be a countable submodel. If $a$ is wide over $M$ and $b$ is random over $M\cup \{a\}$ then $a$ is wide over $M\cup \{b\}$. 
\end{fact}

\subsection{Connected components and the stabilizer theorem}

We recall that a set is type-definable over $B$ if it is an intersection of definable sets over $B$.

We denote with $G_M^{00}$ the intersection of all type-definable subgroups (over $M$) of bounded index, so $G_M^{00}$ is the smallest such. In \cite{hrushovski2012stable}, Hrushovski proved that $G_M^{00}$ is highly related to the notion of measure-teoretical stabilizer.

\begin{defn}
    Given a type $p$ over $M$, we define $$St(p)=\{g\in G \mid p\cap g\cdot p \text{ is wide  }\}$$ where $g\cdot p$ is the partial type over $M,g$ obtained by translating definable sets in $p$ by $g$, that is, the formula $\varphi(x) $ belongs to $ p$ if and only if $\varphi(g^{-1}x) $ belongs to $g\cdot p$.

    We write $Stab(p)$ for the subgroup generated by $St(p)$. 
\end{defn}

For example, if $p$ is not wide, then $St(p)$ is empty; if $p$ is wide, then $1$ is trivially in $St(p)$.

\begin{thm}[{\cite[Theorem 3.5]{hrushovski2012stable}} \& {\cite[Theorem 2.12]{montenegro2020stabilizers}}]\label{T stabilizer}
     If $p$ is a wide type over $M$, then    $$G_M^{00}=Stab(p)=St(p)^2=(p\cdot p^{-1})^2.$$
    Moreover, the subgroup $G_M^{00}$ is wide and every $g$ in $G_M^{00}$ wide over $M$ is contained in $St(p)$.
\end{thm}

By definition of $St(p)$ we deduce that if $g$ is a wide element of $G_M^{00}$, then $p\cap g\cdot p$ is wide. So if $A$ is $M$-definable and it has positive measure, then it can be completed to a wide type $p$; applying the stabilizer theorem to $p$ we get the following.

\begin{coroll}\label{coroll A cap A+g}
    Let $A$ be an $M$-definable set of positive measure and $g$ in $G^{00}_M$ wide over $M$. The set $A\cap g\cdot A$ has positive measure.
\end{coroll}

In \cite{martin2024complete} Martin-Pizarro and Palacín further explored the use of random elements in model-theoretical extensions of amenable groups in order to find additive configurations. We will make a large use of the following result.

\begin{fact}[{\cite[Corollary 3.2]{martin2024complete}}]\label{fact stesso coset}
    If $h$ and $g$ are two elements wide over $M$ in the same coset of $G_M^{00}$ and $A$ is $M$-definable, then $\mu(A \cap h\cdot A)=\mu(A\cap g\cdot A)$.
\end{fact}

\section{Solving equations in pseudo-finite groups}

In this section, we aim to prove Rado's theorem in a more general measure-theoretic context. All groups will be abelian, and we will indicate their operation additively.

\begin{comment}
$\mu$ is an $M$-definable measure, invariant by translations and inverse. We also assume the group $G$ to be abelian (so we indicate the operation with the symbol $+$) and the positivity of the measure to be preserved under scalar multiplication: $$\mu(A)>0 \iff \mu(nA)>0$$ for every non-zero $n$, where $nA=\{na \mid a\in A\}$.

\begin{fact}[{\cite[Satz IV]{Rado}}]
    
    Given a linear polynomial $L(x_1,\dots, x_k)= n_1x_1 +\dots + n_kx_k$ with zero constant term, the following are equivalent.
    \begin{enumerate}
        \item Given any finite coloring of $\N$, there exists a monochromatic solution to $L$;
        \item $\sum_{i\in I}n_i=0$ for some nonempty set $I\subset \{1,\dots, k\}$.
    \end{enumerate}
\end{fact}

We say that a linear equation satisfies Rado's condition if it satisfies the second condition above.
\end{comment}

We will prove that if $X$ is a set of positive measure in $G$ and the linear polynomial $L$ satisfies Rado's condition, then $L$ has solutions in a translate of $X$. The case in which $I$ equals the whole set $\{1,\dots, k\}$ was already proved in (\cite[ Remark 3.15]{martin2024complete}).

\begin{fact}\label{Rado - caso somma 0}
    Let $X$ be an $M$-definable set of positive measure, and $n_1,\dots,n_{s-1}\text{ and } n_{s}$ be integers such that $\sum_{i=1}^{s-1}n_i=n_s$. The set 
    $$\{(x_1,\dots, x_{s-1})\in G^{s-1} \mid x_i\in X \text{ and } n_1x_1 + \dots + n_{s-1} x_{s-1} \in n_sX\}$$ has positive measure.
\end{fact}

%For the sake of completion, we report here the proof. \todo{Insert or not?}

%\begin{proof}
    
%\end{proof}

\begin{defn}
    Given an equation $E(x_1,\dots, x_s)$, we say that the tuple $(a_1,\dots, a_s)$ is a generic solution to $E$ if (up to permutations) $a_{i+1}$ is wide over $a_1,\dots,a_{i-1} $ and $ a_i$ for $i\leq s-2$.
\end{defn}

\begin{rmrk}\label{3:rmrk:PreserveRandomness}
    Assume that for all definable sets $Y$ the following holds: 
    $$\mu(Y)>0 \iff  \mu(N\cdot Y)>0$$ where $N$ is an integer. Then the function sending $a$ to $N\cdot a$ preserve wideness, that is: the element $a$ is wide, resp. random, over some parameter set $B$ if and only if $N\cdot a$ is, for any non-zero integer $n$.
\end{rmrk}

\begin{proof}
Notice that $a$ is in the definable set $ Y$ if and only if $N\cdot a $ is in $ N\cdot Y$; if the positivity of the measure is preserved under multiplications, then we immediately see that $a$ is wide if and only if $N\cdot a$ is. By regularity of the density this result extends to Borel sets in the space of types, from which we obtain that $a$ is random if and only if $N\cdot a $ is.
\end{proof}

\begin{thm}\label{teo rado su traslati} 
Let $X$ be an $M$-definable subset of $G$ with positive measure and $n_1,\dots , n_{s-1} \text{ and } n_{s}$ integers such that  $r:=\sum_{i=1}^{l-1} n_i = -n_l$ and let $k=\sum_{j=l+1}^s n_j$. Assume that multiplication by $N$ preserves randomness, for $N$ among $\{k, r,n_s\}$.
    
    Then, for some $m$ in $G(M)$, the set
    \begin{equation}
        \{(x_1,\dots, x_{s-1})\in G^{s-1} \mid x_i\in (X+m) \text{ and } n_1x_1 + \dots + n_{s-1} x_{s-1}\in n_{s} (X+m)\}
    \end{equation}
    has positive measure. Equivalently, the linear polynomial $n_1x_1 + \dots + n_s x_s$ admits a generic solution in $(X + m)$. 
\end{thm}

\begin{proof}[Proof of \Cref{teo rado su traslati}]
Due to \cref{Rado - caso somma 0}, it suffices to prove the case where the sum of all coefficients is not zero. We start by proving that the set $$\{(x_1,x_{2}, z)\in G^3 \mid x_i\in (X+z) \text{ and } rx_1 + k x_{2}\in r (X+z)\}$$ has positive measure.

    Let $b_2$ be an element of $G_M^{00}$ wide over $M$, and $a$ in $X$ random over $M\cup \{b_2\}$. It follows that $b_2$ belongs to $G_M^{00}\cap (X + g)$ where $g=b_2 - a$; also, $g$ is random over $M\cup \{b_2\}$, as $a$ is, and by \Cref{referee's lemma} it follows that $b_2$ is wide over $g$. 

     Notice that $kb_2$ is a wide element in $G^{00}_M$, so \Cref{coroll A cap A+g} applied to $rX$ implies that $\mu(rX + kb_2 \cap rX)>0$. By translation-invariance we have that also $\mu(r(X + g) + kb_2 \cap r(X+g))>0$. Thus, there exists $$c=rb_3\in r(X + g) + kb_2 \cap r(X+g)$$ random over $M,g, b_2$. For some $b_1$ in $ X+g$ write $rb_1 + kb_2 =c\in r(X+g)$, so $$(b_1, b_2, g)\in \{(x_1,x_{2}, z)\in G^3 \mid x_i\in (X+z) \text{ and } rx_1 + k x_{2}\in r (X+z)\},$$ hence if we show that the tuple $(b_1, b_2, g)$ is wide over $M$, the conclusion follows. 

    \begin{itemize}
        \item $g$ is wide over $M$;
        \item $b_2$ is wide over $M\cup \{g\}$;
        \item $b_1$ is wide over $M\cup \{g,b_2\}$: indeed $c=rb_1 + kb_2$ is random over $M\cup \{g, b_2\}$, hence also $rb_1$ is, and so is $b_1$.
    \end{itemize}

    By \Cref{F wr tuples} the tuple is wide. 
    
\begin{claim}
    We can assume $g$ in $G(M)$.
\end{claim}

\begin{claimproof}
    The wide pair $(b_1, b_2)$ above is contained in the definable set $$\{(x_1,x_2)\in G^2 \mid x_1,x_2\in (X+g) \text{ and } rx_1 + k x_2\in r (X+ g)\}=\varphi(x_1,x_2;g),$$ with $\varphi(x_1, x_2 ,z)$ a formula with parameters from $M$, which then must have positive measure, say $\alpha>0$. Fixed $\varepsilon=\frac{\alpha}{2}$, \Cref{rem b in Y} implies the existence of an $M$-definable set $Z$ containing $g$ and contained in $$\{z\in G \mid \alpha - \varepsilon\leq\mu(\varphi(x_1,x_2;z))\leq \alpha + \varepsilon\}.$$ Let $m$ be in $Z(M)$. It follows that $\mu(\varphi(x_1,x_2;m))\geq \alpha - \varepsilon=\frac{\alpha}{2}>0$ which means that $$\{(x_1,x_2)\in G^2 \mid x_1,x_2\in (X+m) \text{ and } rx_1 + k x_2\in r (X+ m)\}$$ has positive measure.
\end{claimproof}

So, from now on, we will write $m$ instead of $g$, to highlight that it belongs to $G(M)$.

We have found a generic solution $(b_1, b_2, b_3)$  to $rx_1 + kx_2=rx_3$ in $(X+m)$. In particular $$rb_3\in r(X +m)\cap r(X +m) + kb_2,$$ hence this intersection must have positive measure. Take now realizations $a_{l+1}, \dots,a_{s-1}$ and $\tilde a_{s}$ of $\tp(b_2/M)$ with $a_{i+1}$ random over $a_{l+1}, \dots , a_{i-1}$ and $a_i$. Notice that, if we denote by $Cos(x)$ the coset of $x$ with respect to $G^{00}_M$, we have that $Cos(kb_2)=Cos(\sum_{j=l+1}^{s-1} n_j a_j + n_s \tilde a_s )$.

    Since both $kb_{2}$ and $\sum_{j=l+1}^{s-1} n_j a_j + n_s \tilde a_s $ are wide, applying \Cref{fact stesso coset} we deduce 
    \begin{equation*}
        \mu(r(X +m)\cap [r(X + m ) - \sum_{i=l+1}^{s-1} n_i a_i - n_s \tilde a _s])= \mu(r(X + m)\cap [r(X + m ) - kb_2])>0.
    \end{equation*}
So, there exists an element $r\tilde a$ in $r(X + m)\cap [r(X + m ) - \sum_{i=l+1}^{s-1}  n_i a_i- n_s \tilde a _s]$, random over $M\cup \{a_{l+1},\dots, \tilde a_s\}$ (and notice that by hypothesis also $\tilde a$ is random over the same parameters). 

This means that for some $\tilde a ' $, the tuple $(\tilde a, \tilde a',  a_{l+1}, \dots , \tilde a_{s})$ is a generic solution in $(X+m)$ to the equation $$rx_1 - rx_2 + n_{l+1} x_{l+1} + \dots + n_s x_s=0.$$

Take now realizations $a_1,\dots, a_{l-2}\text{ and } a_{l-1}$ of $\tp(\tilde a /M)$ wide one over the previous ones and $\{a_{l+1}, \dots, \tilde a_s\}$. Then $\sum_{i=1}^{l-1} n_i a_i$ and $r\tilde a$ are both wide over $M$ and belong the same coset of $G^{00}_M$. It follows from \Cref{fact stesso coset} that the sets $$n_s (X +m) \cap [r(X+m) -  \sum_{i=1}^{l-1} n_i a_i- \sum_{j=l+1}^{s-1} n_j a_j ]$$ and $$n_s (X +m) \cap [r(X +m) -  r\tilde a- \sum_{j=l+1}^{s-1} n_j a_j ]$$ have the same measure, positive since the latter set contains $n_s\tilde a_s$. Take thus a random (over everything else) element $n_s a_s$ in that intersection: by hypothesis then also $a_s$ will be a random element, and one sees that for some appropriate $a_l$, the tuple $(a_1,\dots, a_s)$ is a generic solution to the original equation. This concludes the proof.
\end{proof}

\begin{comment}
\begin{rmrk}\label{3:rmrk:PreservareLaMisuraNonServe}\todo{I still believe that it's cleaner to keep the this here, so that the theorem's text is shorter}
    Notice that for the proof we do not need the measure $\mu$ to satisfy $$\mu(A)>0\iff \mu(nA)$$ for all integers $n$: it is enough to have this property for all $n$ in between $-N$ and $N$, for some $N$ which explicitly depends on $L$ (for example, $N$ can be the sum of all coefficients $L$, taken with the positive sign).  
\end{rmrk}
\end{comment}

The above proofs yields immediately the following result. 

\begin{coroll}\label{3:Prop:RadoSenzaTraslati}
    Let $X$ be an $M$-definable subset such that $X\cap G_M^{00}$ is wide (as a partial type). For any integers $n_1,\dots, n_s$ satisfying Rado's condition, the following set has positive measure:
    \begin{equation*}
        \{(x_1,\dots, x_{s-1})\in G^{s-1} \mid x_i\in X \text{ and } n_1x_1 + \dots + n_{s-1}x_{s-1}\in n_s X\}.
    \end{equation*}
\end{coroll}

\begin{comment}
\begin{proof}

    First, the conclusion holds when $\sum_{i=1}^s n_i =0$ because of \cite[Remark 3.15]{martin2024complete}. 

    For the other cases, the proof is just an adaptation of the previous one. Using the notation therein, we can solve $rx_1 + kx_2=rx_2$ generically in $X$ (the translation was only needed to intersect the connected component). Then the proof works without any changes.
\end{proof}
\end{comment}

\subsection{Some consequences for finite and pseudo-finite groups}

Given a family of sets $\{M_n\}_{n\in \N}$ and a non-principal ultrafilter $\mathcal U$ on $\N$, we indicate with $\prod_\mathcal{U} M_n$ the ultraproduct of the family. If each $M_n$ is an $\mathcal L$ structure, then $\prod_\mathcal{U} M_n$ is an $\mathcal L$ structure as well due to Łoś' theorem (see \cite[Theorem 2.1.1]{macpherson2018model}), and it is $\aleph_1$ saturated. We call \textit{internal} a set of the form $X=\prod_\mathcal{U} X_n$ with each $X_n$ a subset of $M_n$. If each $M_n$ is finite and $X=\prod_\mathcal{U}X_n$ is an internal set, we define its cardinality as the class of the sequence $(|X_n|)_n$ in $^*\N=\prod_\mathcal{U} \N$. 

\begin{defn}\label{3:def:Loeb}
    Let $M=\prod_\mathcal{U} M_n$ be an ultraproduct of finite structures, and $X=\prod_\mathcal{U} X_n$ be an internal set. We define its Loeb measure $$\mu(X)=st(\frac{|X|}{|M|})$$ where $st:{}^*\R=\prod_\mathcal{U}\R \to \R$ is the standard part function.
\end{defn}

The Loeb measure is a $\sigma$-additive probability measure defined over the $\sigma$-algebra generated by internal sets. Similarly, we can build Loeb measures for the Cartesian powers of $M$: in dimension $n$, it will be given by $$\mu_n(X)=st(\frac{|X|}{|M^n|}).$$

As shown in \cite[Theorem 2.10.3]{tao2014hilbert}, these measures verify Fubini.

\begin{defn}
    A pseudo-finite group is an infinite group $G$ isomorphic to an ultraproduct $\prod_\mathcal{U} G_n$ where each $G_n$ is a finite group.
\end{defn}

It is well known that the Loeb measure over pseudo-finite groups is (left and right) translation-invariant; moreover, following \cite[Section 2.6]{hrushovski2012stable} we can choose a countable language $\mathcal L$ such that each $\mu_n$ is $\emptyset$-definable. Therefore, we can apply the results from this section to pseudo-finite groups, and obtain from that asymptotic consequences for finite groups. 

\begin{coroll}\label{3:co:RadoPerPseudo-finite}
    Let $G$ be a pseudo-finite group and $X$ an internal subset of positive Loeb measure. Consider a linear polynomial $L=n_1x_1+ \dots + n_s x_s$ with $-n_l=\sum_{i=1}^{l-1} n_i$ and assume that for every integer $N$ among $\{n_l, n_s, \sum_{j=l+1}^s n _j\}$ and for every internal set $Y$, in $G$ it holds $$\mu(Y)>0 \iff \mu(N\cdot Y).$$
    Then the set of solutions to $L$ \begin{equation*}
        \{(x_1,\dots, x_{s-1})\in G^{s-1} \mid x_i\in (X+m) \text{ and } n_1x_1 + \dots + n_{s-1} x_{s-1}\in n_{s} (X+m)\}
    \end{equation*} has positive $(s-1)$-dimensional Loeb measure.
\end{coroll}

\begin{coroll}\label{3:co:RadoPerFiniti}
    Fix a natural $r$ and a real $\epsilon$, and let $L(x_1,\dots, x_s)=n_1x_1 + \dots + n_s x_s$ be a linear polynomial over $\Z$ with $-n_l=\sum_{i=1}^{l-1} n_i$. There exist a $c$ and a $\delta$ such that the following holds. 
    
    Let $G$ be a finite abelian group of cardinality  greater then $c$ such that, for every integer $N$ among $\{n_l, n_s, \sum_{j=l+1}^s n _j\}$ the map $x \mapsto Nx$ has fibers of cardinality bounded by $r$, and let $X$ be a subset of $G$ of cardinality at least $\epsilon c$. Then for some $g$ in $G$, the set \begin{equation*}
        \{(x_1,\dots, x_{s-1})\in G^{s-1} \mid x_i\in (X+g) \text{ and } n_1x_1 + \dots + n_{s-1} x_{s-1}\in n_s (X+g)\}
    \end{equation*}
    has cardinality at least $\delta c^{s-1}$.
\end{coroll}

For example, if $|G|$ is coprime with the product $n_l\cdot n_s\cdot\sum_{j=l+1}^s n_j $, then multiplication by $N$ is bijective, and the theorem applies.

\begin{proof}
    Suppose this is not the case. Then there exist $r, \epsilon $ and a linear polynomial $L$ such that for every $c$ and $\delta$ there are a finite abelian group $G_{c, \delta}$ and a subset $X_{c, \delta}$ with the following properties: \begin{itemize}
        \item $|G_{c, \delta}|\geq c$;
        \item $|X_{c, \delta}|\geq \epsilon c$;
        \item the maps $x\mapsto Nx$ have fibers bounded by $r$;
        \item for all $g$ in $G_{c,\delta}$, the set \begin{equation*}
        S_{c,\delta, g}=\{(x_1,\dots, x_{s-1})\in G_{c, \delta}^{s-1} \mid x_i\in (X_{c, \delta} - g) \text{ and } n_1x_1 + \dots + n_{s-1} x_{s-1}\in n_s (X_{c, \delta}+g)\}
    \end{equation*} has cardinality smaller than $\delta c^{s-1}$.
    \end{itemize}

    Fix a non-principal ultrafilter $\mathcal U$ on $\N$, and let $G$ and $X$ be the ultraproducts $$G=\prod_\mathcal{U} G_{c, \frac{1}{c}}\ ,\ X=\prod_{\mathcal U} X_{c, \frac{1}{c}}.$$
    So $G$ is a pseudo-finite group and $X$ is an internal set of Loeb measure at least $\epsilon$. 

    Now, the fact that $x\mapsto Nx$ have fibers bounded by $r$ for $N$ in $\{n_l, n_s, \sum_{j=l+1}^s n_j\}$ implies that, for all internal sets $B$, $$\mu(B)>0 \iff \mu(N\cdot B)>0.$$ So we can apply \Cref{3:co:RadoPerPseudo-finite}, and obtain that for some $g=[(g_n)_n]$, the set \begin{equation*}
        S=\{(x_1,\dots, x_{s-1})\in G^{s-1} \mid x_i\in (X - g) \text{ and } n_1x_1 + \dots + n_{s-1} x_{s-1}\in n_s (X+g)\}
    \end{equation*} has positive $(s-1)$ Loeb measure, say larger than $\eta$. But notice that $$S=\prod_\mathcal{U} S_{c, \frac{1}{c}, g_c}$$ and this would imply that $$\{c \in \mathbb N\mid |S_{c, \frac{1}{c}, g_c} |> \eta c^{s-1}\}\in \mathcal U.$$
    
    Since by construction $|S_{c, \frac{1}{c}, g_c}|<\frac{1}{c}\cdot c^{s-1}$, we must have $\eta < \frac{1}{c}$. In particular, the set above is finite, as it is contained in $[1, \frac{1}{\eta}]\cap \N$, and so it cannot belong to any non-principal ultrafilter.     
\end{proof}

Our method also proves the finitistic version of (the sufficiency part of) Rado's theorem for partition regularity of linear equations (cf. \cite{RadoFinite}, Theorem 3.4).

\begin{thm}
Let $L(x_1,\dots ,x_s)=n_1x_1 + \dots + n_s x_s$ be a linear polynomial with $n_l =-\sum_{i=1}^{l-1} n_i$. For each $m$ and $r$ there exists a $c$ such that the following holds. 
    
    For every finite abelian group $G$ with the properties:
    \begin{itemize}
        \item $|G|\geq c$;
        \item for every $N$ in $\{n_l, n_s, \sum_{l+1}^s n_j\}$, the map $x \mapsto Nx$ has fibers of cardinality bounded by $r$,
    \end{itemize}
    and for every $m$-coloring of $G$ there exists a monochromatic solution to $L=0$. 
\end{thm}

\begin{proof}
    Assume this is not the case: hence, there exist an $m$ and $r$ such that for every $c$ we have a group $G_c$ of cardinality at least $c$, with $x\mapsto Nx$ of small fibers, and an $m$-coloring of $G$ with no monochromatic solutions to $L$. 

    Consider the ultraproduct $G=\prod_{\mathcal U}G_c$ and notice how we can define an $m$-coloring of $G$ just taking the ultralimits of the $m$-colorings of all $G_c$'s; by construction, $L$ will not have any monochromatic solution with respect to this coloring. Fix a countable elementary substructure $M$ and let $g$ be a wide (over $M$) element in $G^{00}_M$, and suppose $C$ is the color containing $g$. As we observed in \Cref{3:Prop:RadoSenzaTraslati}, we can find a tuple $a$ in $C$ (therefore monochromatic) such that $L(a)=0$, against our assumption. 
\end{proof}

\section{Solving equations in pseudo-finite fields}

We now pass to the study of pseudo-finite fields. These are ultraproducts of finite fields; as for the group case, their Cartesian products can be equipped with the Loeb measure, which will be translation invariant $\sigma$-additive measures on internal sets that satisfies Fubini. Moreover, we also have invariance under multiplication: for every internal set $E$ and for every nonzero element $x$, the set $x\cdot E$ has the same measure of $E$. Since our arguments work in a more general setting, we collect our assumptions below.

\begin{ass}
    When not specified, $F$ will be a field definable in an enough saturated model $\U$ of a complete, countable first-order theory, with $M$ as a countable model. We will assume each Cartesian power of $F$ is equipped with a definable Keisler measure $\mu_n$ invariant under translations, multiplications and inverse, and they satisfy Fubini. We also assume that all the multiplicative subgroups $F^{\times n}=\{g^n\mid g\in F\}$ have bounded (equiv. finite) index. 
\end{ass}

\begin{notation}
    As usual, we will indicate with $F^*$ the multiplicative group of invertible elements of $F$, and we will write $F^+$ for intending $F$ as an additive group. For example, $(F^+)^{00}_M$ will indicate the additive connected component over $M$, while $(F^*)^{00}_M$ will indicate the multiplicative one. When writing stabilizers, if needed, we will highlight the operation writing $Stab^+(p)$ for the additive stabilizer, and $Stab^*(p)$ for the multiplicative one.
    
    We also write $A^{\times k}$ to indicate the set $\{a^k \mid a\in A\}$, to distinguish it from the Cartesian power $A^k=A\times \dots \times A$. 
\end{notation}

We start by showing a variation of \cite[Lemma 3.2]{pillay1998supersimple}, in presence of a definable measure. The proof is a verbatim adaptation of the original one.

\begin{lemma}\label{4:le:TcapF+}
    If $T$ is an $\emptyset$-definable multiplicative subgroup of finite index, then $T\cap (F^+)^{00}_M$ is wide over $M$.
\end{lemma}

\begin{proof}
     Since $T$ is definable and of finite index and $M$ is an elementary substructure, each coset has a representative in $F(M)$.     
Note that, if $c$ is an element of $F(M)$, then multiplication by $c$ is a definable bijection that permutes definable subgroups, hence $c\cdot (F^+)^{00}_M= (F^+)^{00}_M$. 

Let $a$ be an element in $(F^+)^{00}_M$ wide over $M$ that lies in a coset $c\cdot T$ with $c$ in $F(M)$. The element $c^{-1}\cdot a$ is wide over $M$ and it is contained in $T \cap c^{-1}\cdot(F^+)_M^{00}=T \cap(F^+)_M^{00}$.
\end{proof}

This lemma can be generalized to the corollary below. Since we will not need it for our purposes, we omit its proof (one could adapt \cite[Lemma 2.4]{martin2004elliptic}). 

\begin{coroll}\label{4:thm:IntFraF+eF*}
    Let $A$ be a small subset of $\U$. Any coset of $(F^+)^{00}_A$ intersects any coset of $(F^*)^{00}_A$ in a point that is wide over $A$. 
\end{coroll}

We can now solve nonlinear equations in large sets.

\begin{comment}
\begin{thm}
    Let $A$ be an $M$-definable subset of $F$ of positive measure. Then there exist $x$ and $y$ in $A$ such that $x-y$ lies in $F^2=\{e^2 \mid e\in F\}$. 
\end{thm}

\begin{rmrk}
    From the proof it will be evident that there is not just a single pair $(x,y)$, but a whole positive-measure set of such pairs. 
\end{rmrk}

\begin{proof}
    Since $A$ has positive measure, there exists a wide type $p$ in $S(M)$ that contains $A$. Since the squares have finite index in $F$, they will contain the multiplicative connected component. From \Cref{4:thm:IntFraF+eF*} we know that there exists an element $g$, wide over $M$, that lies in the intersection $(F^+)^{00}_M\cap (F^*)^{00}_M$. In particular, $g$ is a square and it belongs to $St^+(p)$ because of the stabilizer theorem. 

    This means that we can find a realization $a_1$ of $p$, wide over $M,g$, such that $a_2=a_1 + g$ also realizes $p$. Since $p$ extends $A$, we have that $a_1$ and $a_2$ are our solutions. \todo{Maybe we can do better using random types. But look also at Daniel's proof of a lot of types that commute.}
\end{proof}

\begin{rmrk}
    In the same way, we can prove the following. Given a set $A$ of positive measure, there exist $x$ and $y$ in $A$ such that their difference is an $n$-th power in $F$.
\end{rmrk}

\end{comment}

\begin{thm}\label{4:thm:x^n-y^n=z^m}
    Assume that the maps $x\mapsto x^n$ and $x\mapsto x^m$ preserve randomness.
    If $X$ is an $M$-definable set of positive measure, then the set $$\{(x, y)\in X^2\mid x^n - y^n\in F^{\times m}\}$$ has positive $\mu_2$ measure.     
\end{thm}

\begin{proof}
    We start noticing that, by our assumptions, the multiplicative subgroup $F^{\times m}$ has finite index, so it contains $(F^*)^{00}_M$. By \Cref{4:le:TcapF+} we deduce that $F^{\times m}$ and $(F^+)^{00}_M$ intersect in a wide element $g^m$. Since $X^{\times n}$ has positive measure (because $x\mapsto x^n $ preserve randomness), it can be completed to a wide type $p$ over $M$. \Cref{T stabilizer} implies that $ p \cap (p + g^m)$ is wide and in particular $$\mu(X^{\times n} \cap (X^{\times n} + g^m))>0.$$ Take a random (over $M\cup\{g\}$) element $a_1^n$ in that intersection: by construction, $a^n - g^m$ lies in $X^{\times n}$ and it is random over $M\cup\{g\}$ as well, hence it has the form $a_2^n $ for some $a_2$ wide over $M\cup\{g\}$. We claim that the pair $(a_1, a_2)$ is wide, which concludes the proof. 

    But $a_1^n$, and so $a_1$, is in particular random over $M$; and we need to prove that $a_2$ is wide over $M\cup \{a_1\}$. This is equivalent to prove that $a_2^n$ is wide over $M\cup \{a_1\}$, which reduces to prove that the difference $a_1^n - a_2^n =g^m$ is wide over $M\cup \{a_1\}$. Since $g$ is wide and $a_1$ random over $M\cup\{g\}$, we deduce from \Cref{referee's lemma} that $g$ (and so $g^m$) is wide over $M\cup \{ a_1\}$. 
\end{proof}

We leave open the question whether one could obtain the same result but with $F^{\times m}$ replaced by $P(F)$, where $P$ is a polynomial in $\Z[x]$ with no constant term. Over the integers and for $n=1$, this was proved by Bergelson and Leibman in \cite[Theorem B]{bergelson1996polynomial}.

\begin{question}
    Let $P(x)$ be a polynomial in $\Z[x]$ with no constant term, and let $X$ be a definable set of positive measure. Does the set $$\{(x, y)\in X^2 \mid x - y \in P(F)\}$$ have positive measure? 
\end{question}

However, similarly to what we did for proving the analog of Rado's theorem in the previous section, we can extend \Cref{4:thm:x^n-y^n=z^m} in two ways: we can solve more complicated equations, and up to translations we can find all elements of the solution in $X$. Putting everything together, we obtain the following result. For simplicity, we will assume that a set $X$ is wide if and only if $X^{\times n}$ is wide, for all $n$ in $\N$: equivalently, that $a$ is wide (resp. random) if and only if $a^n$ is.

\begin{coroll}\label{4:co:EquazioniInGenerale}
    Let $X$ be a wide $M$-definable set, and let $c_1, \dots, c_l, d_1,\dots , d_{s-1}$ and $d_{s} $ be elements of $F$ such that $\sum_{i=1}^l c_i=0$. Then for some $g$ in $F(M)$, the set \begin{equation*}\label{4:eq:InsiemeDiSoluzioni}
        \{(x_1, \dots , x_l, y_1,\dots, y_{s-1} )\in (X+ g)^{l +s -1} \mid \sum_{i=1}^l c_i\cdot x_i^n + \sum_{j=1}^{s-1} d_j\cdot y_j^m\in d_{s} (X+g)^{\times m}\}\end{equation*}
    has positive $\mu_{k+l}$ measure.
\end{coroll}

\begin{proof}
We may assume that all $c_i$'s and $d_j$'s are contained in $F(M)$. 

Notice that we can solve generically $c x^n - cy^n=dz^m $ in a translate $(X + g)$ of $X$. Indeed, it suffices to take $d h^k$ in $(F^+)^{00}_M\cap d\cdot  F^{\times k}$ wide (it exists because of \Cref{4:le:TcapF+}) and choose $g$, random over $h$, such that $h$ is contained in $X + g$. Passing eventually to a larger model $M'$ and applying the stabilizer theorem to $\tp(ch^n/M')$ we get our generic solution. The same argument used in \Cref{teo rado su traslati} proves that such a $g$ can be 
chosen in $F(M)$. 

At this point, it suffices to repeatedly use \Cref{fact stesso coset} as we did in the proof of \Cref{teo rado su traslati} to conclude the argument.
\end{proof}

As in the previous section, we can use these results to deduce statements for finite fields. 

\begin{rmrk}
    In the pseudo-finite context, a set $X$ has positive measure if and only if $X^{\times n}$ has, for every natural $n$. Indeed, this follows from the estimate $$\frac{|X|}{n}\leq |X^{\times n}|\leq |X|.$$ In particular, all power set $F^n$ have finite index bounded by $n$, and \Cref{4:thm:x^n-y^n=z^m} holds for all $X$ of positive measure. 
\end{rmrk}

The following is a generalization of \cite[Theorem B]{mccutcheon2003sarkozy}.

\begin{coroll}
    Let $m $ and $n$ be two positive integers, and $\epsilon$ a real number in $[0,1]$. There exist a natural number $c$ and a real $\delta$ such that the following holds.

    For every finite field $K$ of cardinality greater than $c$ and for every subset $X$ of $K$ with density at least $\epsilon$, the set $$\{x_1, x_2 \in X \mid x_1^n - x_2 ^n \in K^{\times m}\}$$ has density at least $\delta$. 
\end{coroll}

\begin{proof}
    Suppose the statement is not true. Then there exist $m,n$ and $\epsilon$ such that, for every $c$ and $\delta$, there exist a finite field $K_{c, \delta}$ of cardinality at least $c$ and a subset $X_{c,\delta}$ with density at least $\epsilon$ such that $$S_{c, \delta}=\{x_1,x_2\in X_{c, \delta}\mid x_1^n - x_2^n \in K_{c, \delta}^{\times m}\}$$ has density less than $\delta$. In particular, this holds for $\delta=\frac{1}{c}$. 

    Given a non-principal ultrafilter $\mathcal U$ over $\N$, let $K$, $X$ and $S$ be the ultraproducts $$K=\prod_{\mathcal U}K_{c, \frac{1}{c}}\ ,\ X=\prod_{\mathcal U}X_{c, \frac{1}{c}}\ , \ S=\prod_{\mathcal U}S_{c, \frac{1}{c}}.$$

    It easily follows from Łoś' theorem that $X$ is an internal set of $K$ of Loeb measure at least $\epsilon$, that $S$ has 2-dimensional Loeb measure $0$, and that $$S=\{x_1, x_2\in X \mid x_1^n - x_2^n \in K^{\times m}\}. $$ 
    This contradicts \Cref{4:thm:x^n-y^n=z^m} 
\end{proof}

Applying the same argument to \Cref{4:co:EquazioniInGenerale} we obtain the following.

\begin{coroll}\label{4:co:NonLinearSuiFiniti}
    Let $m $ and $n$ be two positive integers, and $\epsilon$ a real number in $[0,1]$. There exist a natural number $c$ and a real $\delta$ such that the following holds.

    For every finite field $K$ of cardinality $\geq c$, for every subset $A$ of $K$ with density at least $\epsilon$, and for every elements $c_1,\dots, c_k, d_1,\dots, d_{l+1}$, the set \begin{equation}
        \{(x_1, \dots , x_k, y_1,\dots, y_l )\in (X + g)^{k+l} \mid \sum_{i=1}^k c_i\cdot x_i^n + \sum_{j=1}^l d_j\cdot y_j^m\in d_{l+1} (X +g)^{\times m}\}\end{equation} has density at least $\delta$. 
\end{coroll}

\nocite{*}
\printbibliography

\end{document}